
\hsize=120mm
\vsize=185mm
\parindent=8mm
\frenchspacing

\def\oh{{\cal O}}
\def\Bbb#1{{\bf #1}}
\def\circle{{\circ}}
\def\pdeg{{p^\circle}}
\def\Gal{\hbox{\rm Gal}}

\def\sqr#1#2{{\vcenter{\vbox{\hrule height.#2pt
   \hbox{\vrule width.#2pt height#1pt \kern#1pt
     \vrule width.#2pt}
   \hrule height.#2pt}}}}
\def\whis{\mathchoice\sqr34\sqr34\sqr{2.1}3\sqr{1.5}3}
\def\qed{\ifmmode \eqno{\whis}\else
{{\unskip\nobreak\hfil\penalty50\hskip 2em\hbox{}\nobreak\hfil
${\whis}$\parfillskip=0pt\par\medskip}}\fi}

\def\pra{\par}
\def\proof{\medskip\noindent{\bf Proof.}\quad}
\def\remark{\medskip\noindent{\bf Remark.}\quad}
\def\notation{\medskip\noindent{\bf Notation.}\quad}
\def\definition{\medskip\noindent{\bf Definition.}\quad}
\def\pth{{\zeta_p}}
\def\myfn{\psi}
\def\ell{l}

\def\bur{1}
\def\bec{2}
\def\con{3}
\def\ecm{4}
\def\gmx{5}
\def\tho{6}
\def\thp{7}
\def\thq{8}
\def\yag{9}

\centerline{\bf The Yagita invariant of general linear groups}
\bigskip
\centerline{H.~H.~Glover, I.~J.~Leary and C.~B.~Thomas}
\bigskip
\noindent
We give a definition of the Yagita invariant $\pdeg(G)$ of an
arbitrary group $G$, and compute $\pdeg(GL_n(\oh))$ for each prime
$p$, where $\oh$ is any integrally closed subring of the complex
numbers $\Bbb C$  (i.e., $\oh$ is
integrally closed in its field of fractions $F$, which is also a
subring of $\Bbb C$).  We also show that $\pdeg(SL_n(\oh))$ is equal
to $\pdeg(GL_n(\oh))$ for $n\geq 2$ except that possibly
$\pdeg(SL_n(\oh))$ is $1/2\pdeg(GL_n(\oh))$ for some small $n$ and
\lq small' rings $\oh$.  Our definition of $\pdeg(G)$ extends both
Yagita's original definition for $G$ finite [\yag] and the
definition given by one of us for $G$ of finite virtual cohomological
dimension (or vcd) [\thq].  For $G$ of finite vcd such that the
Tate-Farrell cohomology $\hat H^*(G)$ of $G$ is $p$-periodic,
$\pdeg(G)$ is equal to the $p$-period.  Hence our results may be
viewed as a generalisation of those of B\"urgisser and Eckmann 
[\bur,\bec], who compute the $p$-periods of $GL_n(\oh)$ and
$SL_n(\oh)$ for various $\oh$ such that these groups have finite vcd
and for all $n$ such that these groups are $p$-periodic.  The methods
we use are similar to those used in [\bur].

\notation Throughout this paper $p$ shall be a fixed prime number, $\oh$ an
integrally closed subring of $\Bbb C$, and $F$ the field of fractions
of $\oh$.  $\pth$ shall be a primitive $p$th root of unity in $\Bbb
C$, and $l$ shall be the degree of $F[\pth]$ as an extension of
$F$.  In [\bur,\bec,\ecm], the notation $\phi_F(p)$ is 
used for $l$.  

\definition Let $G$ be a group.  
If $C$ is an order $p$ subgroup of $G$, then 
$$H^*(BC;\Bbb Z)=\Bbb Z[x]/(px)$$
for $x$ a generator of $H^2$.  Define $n(C)$ (either a positive integer
or infinity) to be the supremum of the integers $n$ such
that the image of $H^*(BG)$ in $H^*(BC)$ is contained in the
subalgebra generated by $x^n$.  Now define the Yagita invariant
$\pdeg(G)$ to be 
$$\pdeg(G)=\hbox{\rm l.c.m.}\{2n(C) : C\leq G,\quad |C|=p\},$$
where the least common multiple of the empty set is equal to 1, 
and the least common multiple of an unbounded set of integers, or of a
set containing infinity, is infinity.  The invariant $\pdeg(G)$ depends
on the prime $p$ as well as on the group $G$, but we have resisted
the temptation to write, for example, $2^\circ(G)$ for $\pdeg(G)$ in
the case when $p=2$.  

The following result is immediate from the definition.  
\proclaim Proposition~1.  Let $f\colon H\rightarrow G$ be a group
homomorphism whose kernel contains no element of order $p$.  
(For example, $f$ might be injective, or the kernel of $f$ might be
torsion-free.)  Then $\pdeg(H)$ divides $\pdeg(G)$.  
\qed

It is easy to see that $\pdeg(G)$ is finite if $G$ is finite, and 
Proposition~1 applied in the case when $G$ is the quotient of $H$ 
by a torsion-free normal subgroup of finite index implies that
$\pdeg(H)$ is finite if $H$ is a group of finite vcd (see [\gmx] for a
more detailed proof of this fact).  Using Chern classes, we shall 
show in Corollary~7 that if $G$ admits a faithful finite-dimensional
representation over $\Bbb C$, then $\pdeg(G)$ is finite.

Some remarks concerning the origin of the definition 
of $\pdeg(G)$ and its geometric motivation are in order.  
N.~Yagita made the above definition for finite groups [\yag], and 
proved that if a finite group $G$ acts freely on a (finite) product 
$S^{m-1}\times\ldots\times S^{m-1}$ of spheres of the same dimension
$m-1$, with trivial action on their homology, then $\pdeg(G)$ divides
$m$.  This result may be viewed as a generalisation of the theorem 
bounding the dimension of a sphere with a free $G$-action in the case 
when the cohomology of $G$ is periodic.  One of us suggested extending
the definition to groups of finite vcd in [\thq].  The connection with
actions on products of spheres extends to this case too.  If $G$ is a group
of finite vcd acting freely, properly discontinuously and with trivial
action on homology, on a finite 
product $S^{m-1}\times\ldots\times S^{m-1}\times \Bbb R^n$, then
$\pdeg(G)$ divides $m$, by the same argument as used in [\yag] for the
case when $G$ is finite.  For arbitrary groups, the connection with 
the geometry of group actions is less clear.  One reason for extending
the definition to arbitrary groups is that this makes our results
easier to state and no harder to prove.

To motivate our work we start by stating the main theorems and some 
corollaries.  The proofs of the main theorems will occupy the
remainder of the paper.  

\proclaim Theorem~2.  As explained under the title \lq Notation', let
$\oh$ be an integrally closed subring of $\Bbb C$ with field of
fractions $F$, and let $l=|F[\pth]:F|$.    
Define $\myfn(t)$, for $t$ a positive real, to be the greatest integer
power of $p$ less than or equal to~$t$.  The Yagita invariant
$\pdeg(GL_n(\oh))$ is given by the following table.  
$$\pdeg(GL_n(\oh))=\cases{1&for $n<\ell$\cr 
2\hbox{\rm l.c.m.}\{m\,:\, \ell|m,\quad m|(p-1),\quad m\leq n\}&for 
$\ell\leq n\leq p-1$\cr
2(p-1)\myfn(n/\ell)&for $n\geq p$\cr}$$

\remark The above theorem holds for all primes $p$, but the statement
can be simplified slightly for $p=2$.  In this case $\ell=1$,
and $\pdeg(GL_n(\oh))=2\myfn(n)$ for all $n\geq 1$.  
Note that $\pdeg(GL_n(\oh))$ depends only on $F$, so that if 
$R$ is any subring of $\Bbb C$ such that $\oh\leq R\leq F$, then 
$\pdeg(GL_n(R))=\pdeg(GL_n(\oh))$.  

Replacing $GL_n(\oh)$ by $SL_n(\oh)$ does not change the Yagita 
invariant, except for small values of $n$ and \lq small' rings 
$\oh$, when it may be necessary to introduce a factor of $1/2$.  Thus 
one obtains the following.  
\proclaim Theorem~3.  With notation as above, for $n\geq 2$,
$\pdeg(SL_n(\oh))$ is equal to $\pdeg(GL_n(\oh))$ except that 
$\pdeg(SL_n(\oh))$ may be equal to ${1\over 2}\pdeg(GL_n(\oh))$ 
in the following cases.
\pra
\narrower
\noindent
a) $p=2$, $n=2$, and $\oh$ contains no square root of $-1$.  
\pra
\noindent
b) $p$ is odd, $n$ has the form $n=2^rl$ for $2^r$ dividing $(p-1)/l$, 
and $\oh$ contains no $n$th root of $-1$.  

\proclaim Corollary~4.  Define the function $\psi$ as in the statement of 
Theorem~2.  Then for $n\geq 2$ and any prime $p$, the 
Yagita invariant of $GL_n(\Bbb Z)$ is as follows.  
$$\pdeg(GL_n(\Bbb Z))=\cases{1&for $n<p-1$,\cr 
2(p-1)\psi(n/(p-1))&for $n\geq p-1$.}$$
For all $n\geq 2$ and all $p$, 
$\pdeg(SL_n(\Bbb Z))=\pdeg(GL_n(\Bbb Z))$, except
that for $p$ odd, $\pdeg(SL_{p-1}(\Bbb Z))=p-1$, and for $p=2$, 
$\pdeg(SL_2(\Bbb Z))=2$.  

\proof The statement for the general linear groups is a special case
of Theorem~2.  The only cases when $\pdeg(SL_n(\Bbb Z))$ is not
completely determined by Theorem~3 are those when $n=p-1$ and when 
$n=p=2$, but in these cases $SL_n(\Bbb Z)$ is $p$-periodic and 
the $p$-period is determined in [\bur].  \qed

The following statement may be proved in the same way as 
Theorems 2~and~3, and could be considered as a corollary but for 
the weaker conditions imposed on the subring $R$ of $\Bbb C$.

\proclaim Theorem~5.  Let $R$ be any subring of $\Bbb C$
containing $\pth$, a primitive $p$th root of 1, and let $\psi$ be as in the 
statement of Theorem~2.  Then the Yagita invariant of $GL_n(R)$ is
given by the following table.  
$$\pdeg(GL_n(R))=\cases{
2\hbox{\rm l.c.m.}\{m\,:\, m|(p-1),\quad m\leq n\}&for 
$n\leq p-1$\cr
2(p-1)\myfn(n)&for $n\geq p$\cr}$$
For $R$ as above, $\pdeg(SL_n(R))=\pdeg(GL_n(R))$ if 
any of the following conditions is satisfied:
\pra\hangindent=\parindent
a) $n\geq\hbox{\rm max}\{p,3\}$,
\pra\hangindent=\parindent
b) $n\geq 2$, $p$ is odd, and $R$ contains a $(p-1)$st root of $-1$,
\pra\hangindent=\parindent
c) $n=2$, $p=2$, and $R$ contains a square root of $-1$.

Our upper bound for $\pdeg(GL_n(\oh))$ uses the Chern classes of the
natural representation in $GL_n(\Bbb C)$, and relies on the following 
proposition.  

\proclaim Proposition~6.  Let $f(X)$ be a polynomial over the field
$\Bbb F_p$, all of whose roots lie in $\Bbb F_p^\times$.  If there is a
polynomial $g$ and integer $n$ such that $f(X)=g(X^n)$, then $n$ has
the form $mp^q$ for some $m$ dividing $p-1$ and some positive
integer~$q$.  

\proof Let $n=mp^q$, where $p$ does not divide $m$.  It remains to
prove that $m$ divides $p-1$.  Now $g(X^n)=g(X^m)^{p^q}$, so without
loss of generality we may assume that $q=0$.  Now if $g(Y)=0$ has
roots $y_1,\ldots,y_k$, then the roots of $g(X^m)=0$ are the roots of 
$y_i-X^m=0$ for each $i$, and because $p$ does not divide $m$ these
polynomials have no repeated roots.  If $y_i$ is not an element of
$\Bbb F_p^\times$ then $y_i-X^m$ can have no roots in $\Bbb F_p$.  It
has exactly $m$ roots in $\Bbb F_p$ if and only if the 
inverse image of $y_i$
under the map $x\mapsto x^m$ from $\Bbb F_p^\times $ to itself has
order $m$, which can only happen if $m$ divides $p-1$.  \qed

We will use the following Corollary in the proof of Theorem~2.  

\proclaim Corollary~7.  With notation as in 
the statement of Theorem~2, let $G$ be a subgroup of $GL_n(F)$. 
Then the Yagita invariant $\pdeg(G)$ divides the number 
given for $\pdeg(GL_n(\oh))$ in Theorem~2.  

\proof For each order $p$ subgroup $C$ of $G$, we bound the number
$n(C)$ occuring in the definition of $\pdeg(G)$ by considering the
image of the cohomology $H^*(BGL_n(\Bbb C))$ of the Lie group 
$GL_n(\Bbb C)$ in $H^*(BC)$ under the inclusion of 
$C$ in $GL_n(\Bbb C)$.  Recall (or see for example [\thp]) that 
$H^*(BGL_n(\Bbb C))$ is a free polynomial algebra on generators 
$c_1,\ldots,c_n$, the {\it universal Chern classes}, where $c_i$ has
degree $2i$.  For a representation $\rho$ of a group $H$ in 
$GL_n(\Bbb C)$, the element $\rho^*(c_i)$ is usually written
$c_i(\rho)$, and called the $i$th Chern class of the representation
$\rho$.  The sum $c.(\rho)$ of all the $c_i(\rho)$ is known as the
total Chern class of $\rho$, and has the property that $c.(\rho\oplus
\rho')= c.(\rho)c.(\rho')$.  

If $H^*(BC)=\Bbb Z[x]/(px)$, then the total
Chern classes of the $p$ distinct 1-dimensional 
complex representations of $C$
are $1+ix$ for $i\in \Bbb F_p$, where the case $i=0$ corresponds to
the trivial representation.  If we fix a generator of $C$ such that 
the 1-dimensional representation of $C$ sending this generator to
$\pth$ has total Chern class $1+x$, then the representation sending
this same generator to $\pth^i$ has total Chern class $1+ix$.  
For $F$ as in the statement, there are
$1+(p-1)/\ell$ irreducible $F$-representations of $C$, the trivial
representation of dimension one and $(p-1)/\ell$ others of dimension $\ell$.
Over $\Bbb C$ the faithful representations split as a direct sum of 
one copy each of the representations sending a fixed generator of 
$C$ to $\pth^{r{i^j}}$, for some $r$, where $i$ generates a subgroup
of $\Bbb F_p^\times$ of order $\ell$ and $0\leq j< \ell$.  
It follows that the total Chern classes of these representations 
have the form $1-ix^\ell$, where $i$
ranges over the $(p-1)/\ell$ distinct $\ell$th roots of unity in $\Bbb
F_p$.  

Now the ring $H^*(BC)\otimes \Bbb F_p$ is isomorphic to 
the free polynomial ring $\Bbb F_p[x]$, and the total Chern class of
the inclusion of $C$ in $GL_n(\Bbb C)$ is a polynomial $f(x)$ 
of degree less than or equal to $n$ satisfying the
hypotheses of Proposition~6.  Moreover, since this inclusion factors
through $GL_n(F)$, the total Chern class can be viewed as a polynomial
$\tilde f(y)$ in $y=x^\ell$ of degree less than or equal to $n/\ell$ 
satisfying the hypotheses of the Proposition.  
This polynomial is not the trivial
polynomial because the inclusion of $C$ is not the trivial
representation.  The integer $n(C)$ (resp.\ $n(C)/\ell$) is the greatest
integer $m$ such that $f(x)$ (resp.\ $\tilde f(y)$) can be expressed as a
polynomial in $x^m$ (resp.\ $y^m$), and so by Proposition~2 each of
$n(C)$ and $n(C)/\ell$ is a divisor of $p-1$ multiplied by a power of
$p$.  Note also that $n(C)$ is at most $n$.  This implies that for
each $C$, the number $n(C)$ divides the bound given in Theorem~2, and
hence so does the l.c.m.\ of all such $n(C)$.  \qed

\remark The discussion of the total Chern classes of
$F$-representations of the cyclic group of order $p$ could be avoided
by quoting the general bounds on the orders of Chern classes of
$F$-representations given in [\ecm].  Corollary~7 also holds for 
groups $G$ having a representation in $GL_n(F)$ which is faithful on
every order $p$ subgroup, by the same proof.  Corollary~7 may be 
compared with Lemma~1.2 of [\gmx].

Our lower bounds for $\pdeg(GL_n(\oh))$ and $\pdeg(SL_n(\oh))$ 
use various finite subgroups, which we define below.  

\definition For $p$ an odd prime, and $m$ a divisor of $p-1$, let
$G_1(p,m)$ be the split metacyclic group $C_p\colon C_m$, where $C_m$ acts 
faithfully on $C_p$, and let $G_2(p,m)$ be the split metacyclic group
$C_p\colon C_{2m}$, where $C_{2m}$ acts via the faithful action of its
quotient $C_m$.  For $p$ an odd prime let $E(p,1)$ be the non-abelian
group of order $p^3$ and exponent $p$, and let $E(2,1)$ be the
dihedral group of order eight.  For any $p$, and any posititve integer
$m$, let $E(p,m)$ be the central product of $m$ copies of $E(p,1)$, so
that $E(p,m)$ is one of the two extraspecial groups of order
$p^{2m+1}$.  In ATLAS notation [\con], the group $E(p,m)$ is called 
$p_+^{1+2m}$.  Since the symbol $p$ is already overworked in this
paper we shall not use the ATLAS notation here.  

\proclaim Lemma~8.  For $p$ an odd prime the Yagita invariants of the
groups defined above are as follows.  
$$\pdeg(G_1(p,m))=\pdeg(G_2(p,m))=2m,\qquad \pdeg(E(p,m))=2p^m$$
For $p=2$, let $Q_8$ be the quaternion group of order 8. Then 
$$\pdeg(E(2,m))=2^{m+1}\quad\hbox{and}\quad \pdeg(Q_8)=4.$$

\proof The groups $G_1(p,m)$ and $G_2(p,m)$ are $p$-periodic, and
their $p$-periods are shown to be as claimed in Proposition~3.1 and 
Section~3.3 of [\bec].  The Yagita
invariant for the extraspecial groups was computed in Example~B of 
Section~2 of [\yag], which in turn used results of [\tho].  
\qed

\proclaim Lemma~9.  For $p$ an odd prime, and $m$ a divisor of $p-1$,
let $n=m\ell/(m,\ell)$.  The group $G_1(p,m)$ embeds in $GL_n(\oh)$, 
and in $SL_{n+1}(\oh)$.  
If either $\ell/(m,\ell)$ is even or $m$ is odd (and greater than 1)
then $G_1(p,m)$ embeds in $SL_n(\oh)$.  If $m$ is even 
and $\oh$ contains an $m$th
root of $-1$, then $G_2(p,m)$ embeds in $SL_n(\oh)$.  
The group $E(p,m)$ embeds in $SL_{\ell p^m}(\oh)$.  
\pra
For $p=2$ and $m>1$ the group $E(2,m)$ embeds in $SL_{2^m}(\Bbb Z)$.  
The group $E(2,1)$ (the dihedral group of order 8) embeds in 
$GL_2(\Bbb Z)$ and in $SL_3(\Bbb Z)$.  The quaternion group of order 
eight embeds in $SL_2(\Bbb Z[i])$.  

\proof First consider the case when $p$ is odd.  
For $G_1(p,m)$ and $G_2(p,m)$ this is essentially contained 
in [\bec] Section~3.2.  If $m$ divides $m'$, then $G_1(p,m)$ is a
subgroup of $G_1(p,m')$, so without loss of generality we may assume 
that $m$ is divisible by $\ell$ and show that in this case $G_1(p,m)$
embeds in $GL_m(\oh)$.  As an $\oh$-module, $\oh[\pth]$ is free of 
rank $\ell$, which is closed under multiplication by $\pth$ and under 
the action of $\Gal(F[\pth]/F)$, which together generate a group of 
$\oh$-linear automorphisms of $\oh[\pth]$ isomorphic to $G_1(p,\ell)$.
If $m$ is a proper multiple of $\ell$, view $G(p,\ell)$ as a subgroup
of $G(p,m)$, and the induced representation $V$ coming from the above
representation is an $m$-dimensional faithful representation of 
$G(p,m)$ over $\oh$.  

It is easy to check that the determinant of the action of an element
of order $p$ of $G_1(p,m)$ on $V$ is equal to 1, and that $V\otimes F$
restricts to a cyclic subgroup of $G_1(p,m)$ of order $m$ as 
a sum of $\ell/(m,\ell)$ copies of the regular representation (use the
normal basis theorem).  Thus if either $m$ is odd or $\ell/(m,\ell)$
is even, the determinant of an element of order $m$ acting on $V$ 
is 1, and so $G_1(p,m)$ is contained in $SL(V)\cong SL_n(\oh)$. 
In any case, the determinant of the action of $G_1(p,m)$ on $V$ has
image contained in $\{\pm 1\}\subseteq \oh^\times$, so the action of 
$G_1(p,m)$ on $V\oplus\Lambda^n(V)$ gives an embedding from $G_1(p,m)$
into $SL_{n+1}(\oh)$.  

If $m$ is even and $l/(m,l)$ is odd, then $\oh$ contains an $m$th root
of $-1$ if and only if $\oh$ contains an $n$th root of $-1$.  In 
this case, if $\mu$ is an $n$th root of $-1$ in $\oh$,
and $A$ of order $p$ and $B$
of order $m$ with ${\rm det}(B)=-1$ generate a subgroup of $GL_n(\oh)$ 
isomorphic to $G_1(p,m)$, then $A$ and $\mu B$ generate a subgroup of 
$SL_n(\oh)$ isomorphic to $G_2(p,m)$.

Now consider the extraspecial group $E(p,m)$ for $p$ odd.
Note that the centre
$Z$ of $E(p,m)$ is cyclic of order $p$, and that $E(p,m)$ has a subgroup 
of index $p^m$ containing $Z$ as a direct factor.  This subgroup has
an $\ell$-dimensional $\oh$-representation which is faithful on $Z$,
and the corresponding induced $E(p,m)$-module is an $\ell
p^m$-dimensional representation which must be faithful (since any
nontrivial normal subgroup of a $p$-subgroup meets the centre
nontrivially).  Over $\Bbb C$, $E(p,m)$ has $(p-1)$ faithful irreducible
representations, each of dimension $p^m$ (arising as induced modules
in the above way for different choices of 1-dimensional modules for 
$Z$).  Using characters it is easy to see that these representations 
restrict to any non-central subgroup of order $p$ as a sum of copies
of the regular representation and to the centre as a sum of copies of
a single irreducible representation, and hence that they 
have image in $SL_{p^m}(\Bbb C)$
(note that $E(p,m)$ has exponent $p$).  The representations over $\oh$ 
constructed above split over $\Bbb C$ into a sum of some of the 
faithful irreducibles, so have determinant~1.  

For $p=2$, there is a unique faithful irreducible complex
representation of $E(2,m)$, which has dimension $2^m$, 
and an argument similar to the above shows
that it is realisable over $\Bbb Z$.  Using characters one can show
that the restriction of this representation to any cyclic subgroup of 
order four lies in $SL_{2^m}(\Bbb Z)$, and that the restriction to 
any non-central subgroup of order two is isomorphic to a sum of 
$2^{m-1}$ copies of the regular representation, so lies in
$SL_{2^m}(\Bbb Z)$ provided that $m\geq 2$.  

The left action of the quaternion group $Q_8$ on $\Bbb Z[i,j,k]$ commutes
with the right action of $\Bbb Z[i]$, giving a faithful representation
of $Q_8$ which has image in $SL_2(\Bbb Z[i])$.  
\qed

\par\noindent
{\bf Proof of Theorem~2.}\quad  Corollary~7 gives an
upper bound for $\pdeg(GL_n(\oh))$.  
For $n\leq p-1$, for each $m\leq
n$ such that $\ell$ divides $m$ and 
$m$ divides $p-1$, Lemma~9 tells us 
that $G_1(p,m)$ occurs as a subgroup 
of $GL_n(\oh)$, and so by Proposition~1 and 
Lemma~8, $\pdeg(GL_n(\oh))$ is divisible by $\pdeg(G_1(m,p))=2m$.  
This gives the bound for $n\leq p-1$, and shows that for $n\geq p$, 
$\pdeg(GL_n(\oh))$ is divisible by $2(p-1)$.  Now (for any $p$) the 
group $E(p,m)$ is a subgroup of $GL_n(\oh)$ 
for each $n\geq 2\ell p^m$, and
has Yagita invariant $2p^m$.  
This gives the $p$-part of the bound for
$n\geq p$.  \qed

\par\noindent
{\bf Proof of Theorem~3.}\quad This is similar to the proof of
Theorem~2.  In the case when 
$n\leq p-1$, the l.c.m.\ occuring in the
expression given for $\pdeg(GL_n(\oh))$ is clearly equal to the 
following expression. 
$$\hbox{\rm l.c.m.}\{m\, : \,m\leq n,  
m = \ell q^r \,\hbox{for some prime $q$},\, 
q^r| (p-1)/\ell\}$$
In other words, we need only 
consider those $m$ of the form $\ell q^r$
for some prime~$q$ such that $q^r$ 
divides $(p-1)/l$.  If $q$ is an odd
prime and $q$ divides $\ell$ exactly 
$s$ times, then $G_1(p,q^{r+s})$ 
is a subgroup of $SL_n(\oh)$ for $n=\ell q^r$, and has Yagita
invariant $2q^{r+s}$.  If 2 divides $\ell$ exactly $s$ times, and 
$2^r\ell$ divides $p-1$, then $G_1(p,2^{r+s})$ is a subgroup of 
$SL_{n+1}(\oh)$ for $n=\ell 2^r$ and has Yagita 
invariant $2^{1+r+s}$.  
From these examples it already 
follows that $\pdeg(SL_n(\oh))$ is 
divisible by $2(p-1)$ for $n\geq p$, and that for $n\leq p-1$, 
$\pdeg(SL_n(\oh))$ is equal to $\pdeg(GL_n(\oh))$ 
except possibly if 
$n$ is of the form $2^r\ell$ and 
is a factor of $p-1$, when the Yagita
invariant for $SL_n(\oh)$ might 
be half the Yagita invariant for
$GL_n(\oh)$.  If $\oh$ contains a 
$2^{r+s}$th root of ${-1}$, or 
equivalently an $n$th root of $-1$ (where $n=\ell 2^r$), then 
$G_2(p,2^{r+s})$ is a subgroup of 
$SL_n(\oh)$, and has Yagita invariant 
$2^{1+r+s}$, so that in this case too
$\pdeg(SL_n(\oh))=\pdeg(GL_n(\oh))$.  
\qed

\noindent
{\bf Proof of Theorem~5.}\quad The groups $GL_n(R)$ and $SL_n(R)$
are subgroups of $GL_n(\Bbb C)$, so their Yagita invariants are 
bounded above by $\pdeg(GL_n(\Bbb C))$.  By the hypothesis on $R$, 
the cyclic group $C_p$ admits a 
faithful representation in $GL_1(R)$.
As in Lemma~9 one may use induced representations of the groups 
$G_i(p,m)$ and $E(p,m)$ to give lower bounds
equal to the above upper bounds.  We leave the details
as an exercise.  
\qed

\remark The methods that we use also gives some information 
concerning the Yagita invariant of the groups $G(\oh)$ for other
algebraic groups~$G$.  We hope to address this question in a future 
publication.  

\medskip
\noindent
{\bf Acknowledgements.} The work of the first named author was 
partially funded by the Centre de Recerca Matem\`atica, and  
the work of the second named author was 
funded by a DGICYT research fellowship at the 
Centre de Recerca Matem\`atica.

\beginsection References.  

\frenchspacing
\def\book#1/#2/#3/#4/#5/{\item{#1} #2, {\it #3,} #4, {\oldstyle #5}.
\par\smallskip}
\def\paper#1/#2/#3/#4/#5/(#6) #7--#8/{\item{#1} 
#2, #3, {\it #4,} {\bf #5}
({\oldstyle#6}) {\oldstyle #7}--{\oldstyle#8}.\par\smallskip}
\def\prepaper#1/#2/#3/#4/#5/{\item{#1} #2, #3, {\it #4} {#5}.
\par\smallskip}

\paper \bur/B. B\"urgisser/On the $p$-periodicity of arithmetic
subgroups of general linear groups/Comment. Math. Helv./55/%
(1980) 499--509/

\paper \bec/B. B\"urgisser and B. Eckmann/The $p$-periodicity of the
groups
\pra
$GL(n,O_S(K))$ and 
$SL(n,O_S(K))$/Mathematika/31/(1984) 89--97/

\book \con/J. H. Conway {\it et. al.\/}/An ATLAS of finite 
groups/Oxford University Press/1985/

\paper \ecm/B. Eckmann and G. Mislin/Chern classes of group
representations over a number field/Compositio 
Math./44/(1981) 41--65/

\paper \gmx/H. H. Glover, G. Mislin and Y. Xia/On 
the Yagita invariant
of mapping class groups/Topology/33/(1994) 557--574/

\paper \tho/C. B. Thomas/A model for the 
classifying space of an extra
special $p$-group/Mathematika/22/(1975) 182--187/

\book \thp/C. B. Thomas/Characteristic classes and the cohomology 
of finite groups/Cambridge University Press/1986/

\paper \thq/C. B. Thomas/Free actions by $p$-groups on products of
spheres and Yagita's invariant $po(G)$/Lecture Notes in 
Math./1375/(1989) 326--338/

\paper \yag/N. Yagita/On the dimension of 
spheres whose product admits
a free action by a non-abelian group/Quart. J. Math. 
Oxford/36/(1985) 117--127/

\bigskip

\noindent
H. H. Glover, Ohio State University, Columbus, Ohio 43210.

\bigskip
\noindent 
I. J. Leary, Centre de Recerca Matem\`atica, 
Institut d'Estudis Catalans, 
E-08193, \hyphenation{Bella-terra} Bellaterra.

\bigskip
\noindent
C. B. Thomas, DPMMS, University of Cambridge, Cambridge CB2 1SB.

\end